\documentstyle[12pt]{amsart}

\def\NZQ{\mathbb}               
\def\NN{{\NZQ N }}

\def\ZZ{{\NZQ Z}}

\def\operator#1#2{\def#1{\mathop{\kern0pt\fam0#2}\nolimits}} 
\def\opn#1#2{\def#1{\mathop{\kern0pt\fam0#2}\nolimits}}  
\operator\chara{char}
\operator\Tor{Tor}
\operator\reg{reg}

\opn\gr{gr}
\opn\Rees{{\mathcal  R}}
\opn\rank{rank}
\def\pot#1#2{#1[\kern-0.28ex[#2]\kern-0.28ex]}

\let\sect=\cap
\let\dirsum=\oplus

\let\to=\rightarrow
\let\To=\longrightarrow
\def\Implies{\ifmmode\Longrightarrow \else
        \unskip${}\Longrightarrow{}$\ignorespaces\fi}
\def\implies{\ifmmode\Rightarrow \else
        \unskip${}\Rightarrow{}$\ignorespaces\fi}

\let\:=\colon

\newtheorem{theorem}{Theorem}[section]
\newtheorem{lemma}[theorem]{Lemma}
\newtheorem{corollary}[theorem]{Corollary}
\newtheorem{proposition}[theorem]{Proposition}
\theoremstyle{definition}

\newtheorem{example}[theorem]{Example}

\newtheorem{definition}[theorem]{Definition}
\let\epsilon=\varepsilon
\let\phi=\varphi
\let\kappa=\varkappa
\opn\ini{in}
\opn\sat{sat}
\def\mm{{\mathfrak m}}

\opn\height{height}
\opn\length{length}
\newtheorem{question}[theorem]{Question}

\def\cocoa {\mbox{\rm C\kern-.13em o\kern-.07 em C\kern-.13em
o\kern-.15em A}}

\def\addots{\mathinner{\mkern1mu\raise1pt\hbox{.}\mkern2mu\raise4pt\hbox{.}
           \mkern2mu\raise7pt\vbox{\kern7pt\hbox{.}}\mkern1mu}}

\begin{document}
\title{Castelnuovo-Mumford regularity of products of ideals}
\author{Aldo Conca  \and  J\"urgen Herzog }
    \maketitle
\section*{Introduction}

Let $R$ be a polynomial ring over a field, $I\subset R$ a graded ideal
and $M$ a finitely  generated graded $R$-module. The highest degree of a
generator of the product
$IM$ is bounded  above by the sum  of  the highest degree of a generator
of
$M$ and the highest degree of a generator of $I$. One may wonder whether
the same relation holds also for the Castelnuovo-Mumford regularity,
that is, whether
\begin{eqnarray}
\label{equa}
\reg(IM)\leq\reg(M)+\reg(I).
\end{eqnarray} This is not the case in general. There are examples
already with
$M=I$ such that $\reg(I^2)>2\reg(I)$, see Sturmfels \cite{St} and Terai
\cite{T}. On the other hand, Chandler \cite{Ch} and  Geramita,
Gimigliano and  Pitteloud
\cite{GGP} have shown that
$\reg(I^k)\leq k\reg(I)$  holds for ideals with  $\dim R/I\leq 1$. In
general one has that  $\reg(I^k)$ is asymptotically a linear function of
$k$, see
\cite{K, CHT}.  If one takes $I=\mm$ and $M$  any graded $R$-module, then
$\reg(\mm M)\leq \reg(M)+1$ holds. So it is  natural to ask whether
    (\ref{equa}) holds  whenever  $I$ is generated by  a regular
$R$-sequence or at  least by a sequence of linear forms. Unfortunately
this is also not the case, even when $M$ is  a monomial ideal with a
linear resolution
     and $I$ is generated by a subset of the variables, see Example
\ref{not}. The purpose of this note is to describe some cases where
(\ref{equa})  is nonetheless valid.

In Section 1 we recall some generalities about regularity and show in
Section 2 that
     (\ref{equa})  is valid for ideals generated by sequences which are
almost regular with respect to
$M$ and regular with respect to $R$, see \ref{complete}. For example,
any generic sequence of  homogeneous forms of length $\leq \dim R$ has
these properties. We also show the validity of  (\ref{equa}) when the
dimension of $I$ is $\leq 1$. The argument is similar as in the
corresponding result of Chandler.

More surprising is the fact, proved in Section 3 (Theorem 3.1), that any
product of ideals of  linear forms has a linear resolution.  This is
obtained as a consequence of a description of a primary decomposition of
such an ideal, see \ref{primdec}.

In Section 4 we consider ideals with linear quotients,  that is, ideals
which can be generated by a minimal system of generators   whose
successive colon ideals  are generated  by linear forms. Examples of
such ideals are stable, and squarefree stable ideals in the sense of
Eliahou-Kervaire \cite{EK} and Aramova-Herzog-Hibi \cite{AHH}, as well
as polymatroidal ideals, as noted in \cite{HT}. Again it turns out that
the property of having  linear quotients is not preserved under taking
products or powers. However we show in Section 5 that products of
polymatroidal ideals are again polymatroidal, and hence have again
linear quotients. This is also  implied by the fact that discrete
polymatroids are just the integer vectors of an integral  polymatroid
(see
\cite[Theorem 3.4]{HH}) and a theorem on polymatroidal sums
\cite[Theorem 3]{W}.

Let $X$ be a generic Hankel matrix and let $I_t$ be the ideal of the
minors of size $t$ of $X$. It has been shown in \cite{C} that $I_2^k$
has a linear resolution for all $k$. Furthermore, it follows from
results in \cite{B} and \cite{C2} that
$I_t^k$ has a linear resolution for all $k$ and for all $t$.  As an
application of the concept of ideals with linear quotients we show in
the  last  section that any product
$I_{t_1}\cdots I_{t_k}$ of ideals of minors of a generic Hankel  matrix
has a linear resolution.
\medskip

Some of the results of this  paper have been conjectured after explicit
computations performed by using the computer algebra system \cocoa
\cite{CNR}. We would like to thank J.~Abbott, A.~Bigatti and  M.~Caboara
for their help and suggestions in doing these computations.

\section{Generalities} Let $K$ be a field and let $R$ be a polynomial
ring over $K$.   Let
$M=\dirsum_{i\in \ZZ}  M_i$ be a finitely  generated graded
$R$-module. For every $i\in \NN$  one defines
$$t^R_i(M)=\max\{ j \,  | \,    \beta^R_{ij}(M)\neq 0\}$$ where
$\beta^R_{ij}(M)$ is the $ij$th graded Betti number of $M$ as an
$R$-module, i.e.
$$\beta^R_{ij}(M)=\dim_K \Tor^R_i(M,K)_j$$ and $t^R_i(M)=-\infty$ if it
happens that
$\Tor^R_i(M,K)=0$. The {\em Castelnuovo-Mumford regularity}
$\reg(M)$ of $M$ is given by
    $$\reg(M)=\sup \{ t^R_i(M)-i \,  : i\in \NN \} $$ The {\em initial
degree}  of a non-zero graded $R$-module $M$ is the least
$i$ such that $M_i\neq 0$. An  $R$-module $M$ has a {\em linear
resolution} if  its regularity is equal to its initial degree. In other
words, $M$ has linear resolution if  its minimal generators have all the
same degree and the matrices of the minimal free resolution of $M$ over
$R$ have entries of degree $1$.

A short exact sequence
$$0 \to N \to M \to P \to 0$$ of graded $R$-modules yields a long exact
sequence of
$\Tor$-modules
$$ \cdots\to  \Tor_{i+1}^R(P,K)\to \Tor_i^R(N,K) \to \Tor_i^R(M,K)   \to
\Tor_i^R(P,K)\to  \cdots $$ It follows that
\begin{eqnarray}
\label{basic}
\reg(M)&\leq & \max\{\reg(P) ,  \reg(N)\}\nonumber  \\
\reg(N)&\leq &\max\{ \reg(M),\reg(P)+1 \}  \\
\reg(P)&\leq &\max\{\reg(N)-1, \reg(M)\}\nonumber
\end{eqnarray}

Let $N$ be a graded module of finite length. We set $s(N)=\max\{s\:
N_s\neq 0\}$. One has (see \cite[Cor.20.19]{E}):

\begin{lemma}
\label{elem-1} Let $N$ be  a graded $R$-module of finite length. Then:
\begin{itemize}
\item[(a)]   $\reg(N)=s(N)$
\item[(b)]  If $0\to N\to M\to P\to 0$ is a short exact sequence of
graded modules then $\reg(M)= \max\{\reg(P), s(N)\}$.
\end{itemize}
\end{lemma}

Let $M$ be a graded $R$-module. A homogeneous element $x\in R$ of
degree $d$  is called {\em almost regular  on $M$} if the multiplication
map $x\: M_{i-d}\to M_i$ is injective for all
$i\gg 0$. Let
$N=H^0_\mm(M)$, i.e.\  $N=\{ a\in M : \mm^k a=0\quad \text{for
some}\quad  k\}$.  Then
$x$ is almost regular for $M$ if and only if $x$ is a non-zerodivisor on
$M/N$.

A sequence $x_1,\ldots,x_m$ of homogeneous elements of $R$ is called an
{\em almost regular
$M$-sequence} if  $x_i$ is almost regular on $M/(x_1,\ldots,x_{i-1})M$
for
$i=1,\ldots,n$.

\begin{proposition}
\label{elem-2} Let $M$ be a graded $R$-module  and  $x\in R$  an almost
regular element on
$M$ of degree
$d$. Set
$N=H^0_{\mm}(M)$. Then $$\reg(M)=  \max\{ \reg(M/xM)-d+1, s(N)\}. $$
\end{proposition}
\begin{pf}  Set $a=\reg(M), b=\reg(M/xM), c=\reg(xM)$ and $s=s(N)$. We
have to show that  $a=\max\{b-d+1, s \}$.  Let $W=(0:_Mx)$; then
$W\subset N$ and
$s(W)=s(N)=s$. We obtain two exact sequences
\[ 0\to W(-d)\to M(-d)\to xM\to 0,
\] and
\[ 0\to xM\to M\to M/xM\to 0.
\]

   From  \ref{elem-1} and the first exact sequence we have
$$(i)\ \ a=\max\{c-d, s\}.$$     From the second  exact sequence and
(\ref{basic}) we have:
$$   (ii) \ \ c\leq \max\{a,b+1\}, \ \ \ \   (iii) \ \ b\leq
\max\{a,c-1\}.$$ By (i) and (ii) we have $a=\max\{c-d, s\}\leq
\max\{a-d, b+1-d,  s\}$ which implies
$a\leq \max\{b+1-d,  s\}$.  By (iii) and (i) we have $b\leq
\max\{a,c-1\}=\max\{c-d,s,c-1\}=\max\{s,c-1\}$.  Hence $\max\{b+1-d,
s\}\leq
\max\{s+1-d,c-d ,s\}=\max\{c-d ,  s\}=a$.
\end{pf}

Given a homogeneous ideal $I$ in a polynomial ring $R$  and a graded
$R$-module $M$, one defines  the saturation $(IM)^{sat}$ of $IM$ as
follows:

$$(IM)^{sat}=\{ x\in M \:  \mm^kx \subset IM\quad  \text{for some}\quad
k\}$$ and the {\em saturation degree} $\sat(IM)$ the smallest index $j$
such that $IM$ and $(IM)^{sat}$ coincide from degree $j$ on. In other
words,
$\sat(IM)=s((IM)^{sat}/IM)+1$.   Note that
$H^0_\mm(M/IM)=(IM)^{sat}/IM$,  and hence $\sat(IM)$ is the smallest
index $j$ such that $H^0_\mm(M/IM)$ vanishes from degree $j$ on. As a
consequence of
\ref{elem-2} we  have

\begin{corollary}
\label{sat} Let $I\subset R$ be a homogeneous ideal, and let
$x\in R$ be a linear form which is almost regular on $R/I$. Then
$\reg(I)=
\max\{\reg(I+(x)), \sat(I)\}$.
\end{corollary}

\section{Regularity of products of ideals and  modules}

Given a graded $R$-module $M$   and a homogeneous ideal $I\subset R$,
the purpose of this section is to discuss cases in  which the inequality
(\ref{equa})     holds.  We mentioned already in the introduction that
this is not always the case. On the other hand, if one takes $I=\mm$,
where $\mm$ is the graded maximal ideal of $R$,  then
$\reg(\mm M)\leq \reg(M)+1$ and hence (\ref{equa}) holds. So it is
natural to ask whether (\ref{equa}) holds  in case $I$ is generated by an
$R$-regular sequence. Unfortunately this is not true, even when
$I$ is generated by linear forms, as the following example shows.

\begin{example}
\label{not} Let $R=K[a,b,c,d]$,  and let $J=(a^2b,abc,bcd,cd^2)$. The
resolution of
$J$ is
$$0\to R^{3}(-4) \to R^{4}(-3) \to J \to 0$$ It follows that
$\reg(J)=3$. If we take
$I=(b,c)$ then the resolution of $IJ$ is
$$0 \to  R(-8)  \to  R^{3}(-6)\dirsum R^{2}(-7) \to R^{10}(-5)\dirsum
R(-6) \to R^{8}(-4) \to IJ\to 0$$ The  non-linear minimal syzygy among
the generators of $IJ$ is
$a^2(bcd^2)-d^2(a^2bc)=0$. Anyway,  $\reg(IJ)=5$.
\end{example}

On the other hand, one has

\begin{theorem}
\label{almost} Let $R$ be a polynomial ring, and $I\subset R$ be an
ideal which is generated by an almost  regular $M$-sequence
$x_1,\ldots,x_m$ with
$\deg x_i =d_i$ for
$i=1,\ldots,m$. Then
$\reg(IM)\leq \reg(M)+d_1+d_2+\cdots +d_m-m+1$.
\end{theorem}

\begin{pf} Set $N_i=H^0_\mm(M/(x_1,\ldots,x_{i-1})M)$. Then by
\ref{elem-2} we have
\[
\reg(M/(x_1,\ldots,x_{i-1})M)=\max\{\reg(M/(x_1,\ldots,x_i)M)-d_i+1,s(N_i)\}
\] for all $i=1,\ldots,m$. This implies that
\[
\reg(M)=\max\{\reg(M/IM)-d_1-d_2-\cdots-d_m+m, s_1,\ldots,s_m\},
\] where $s_i=s(N_i)-d_1-\cdots -d_{i-1}+(i-1)$.

Thus we see that
\[
\reg(M/IM)\leq \reg(M)+d_1+\cdots +d_m-m.
\] Now
\begin{eqnarray*}
\reg(IM)&\leq &\max\{\reg(M/IM)+1, \reg(M)\}\\ &\leq
&\max\{\reg(M)+d_1+\cdots +d_m-m+1, \reg(M)\}\\ &=&\reg(M)+d_1+\cdots
+d_m-m+1.
\end{eqnarray*}
\end{pf}

\begin{corollary}
\label{complete} Suppose that, in addition to the assumptions of {\em
\ref{almost}},
$x_1,\ldots,x_m$ is a  regular $R$-sequence. Then $\reg(IM)\leq
\reg{M}+\reg(I)$.
\end{corollary}

\begin{pf} For the proof we just note that
$\reg(I)=d_1+d_2+\cdots+d_m-m+1$ if
$x_1,\ldots, x_m$ is a  regular $R$-sequence.
\end{pf}

\begin{corollary}
\label{generic} Let $I$ be an ideal generated by a generic sequence of
homogeneous forms of length
$\leq \dim R$. Then
$\reg(IM)\leq \reg(M)+\reg(I)$.
\end{corollary}

\begin{pf} A generic sequence is an almost regular sequence on $M$ and a
regular sequence on $R$.
\end{pf}

The following result generalizes a theorem of \cite{Ch} and
\cite{GGP}, and is another case in  which the inequality (\ref{equa})
holds

\begin{theorem}
\label{smalldim} Let $R$ be a polynomial ring, and let $I$ be a graded
ideal with $\dim R/I\leq 1$. Then for any finitely generated graded
$R$-module $M$ we have
$$\reg(IM)\leq\reg(M)+\reg(I).$$
\end{theorem}

\begin{pf} The proof follows very much the line of arguments of
\cite{Ch}. 

Let $x\in R_1$ be an element which is almost regular on $M$, $M/IM$ and $R/I$.
We first show that 
\begin{eqnarray} 
\label{satIM}
\sat(IM)\leq \reg(M)+\reg(I).
\end{eqnarray}
We set $r=\reg(M)$ and $t=\reg(I)$. Since $(IM)^{sat}/IM$ and $(IM:_Mx)/IM$ have the same socle, it suffices to
show that if $f\in  M$ is homogeneous of degree $> r+t$  with $xf\in IM$, then $f\in
IM$.

Suppose that $f=\sum_i f_im_i$ and $xf=\sum_ig_im_i$ with $g_i\in I$.
Then
$\sum_i(xf_i-g_i)m_i=0$. Consider the exact sequence
\begin{eqnarray*} 0\To U @>>> F @>\epsilon >> M\To 0
\end{eqnarray*}  where $F$ is free with basis $e_1,\ldots, e_k$ and
$\epsilon(e_i)=m_i$. Then
$\sum_i(xf_i-g_i)e_i\in U$. Let $u_1,\ldots, u_l$ be a homogeneous
system of generators of
$U$, and  $u_j=\sum_ja_{ij}e_i$.  Then $\sum_i(xf_i-g_i)e_i=\sum_j
k_ju_j=\sum_i(\sum_ja_{ij}k_j)e_i$, so that
$xf_i-g_i=\sum_ja_{ij}k_j$. Note that
$\deg  k_j>r+t+1-\deg u_j\geq t$. Hence, $k_j\in I+(x)$, since
$(I+(x))_i=R_i$ for
$i> t$. Thus
$k_j=xp_j+q_j$ with $q_j\in I$. This yields
\[ x(f_i-\sum_ja_{ij}p_j)=g_i+\sum_ja_{ij}q_j
\] This equation implies that $f_i-\sum_ja_{ij}p_j\in I^{sat}$. However,
since
$\sat(I)\leq
\reg(I)=t$  and $\deg(f_i-\sum_ja_{ij}p_j)>t$, it follows that
$f_i-\sum_ja_{ij}p_j\in I$. We  conclude that
$f=\sum_i(f_i-\sum_ja_{ij}p_j)m_i\in IM$. This concludes the proof of (\ref{satIM}).

In order to prove the theorem we assume first that $\dim  M/IM=0$. By ($\ref{basic}$) we
have $\reg(IM)\leq \max\{\reg(M), \reg(M/IM)+1\}$. Hence it  suffices to show that
$\reg(M/IM)\leq \reg(M)+\reg(I)-1$. Since $\reg(M/IM)=s(M/IM)$ by \ref{elem-1}, and since $s(M/IM)=\sat(IM)-1$, this 
follows from (\ref{satIM}).

Now we assume that $\dim M/IM=1$. Set $N =M/xM$. Then  Proposition \ref{elem-2} implies
\begin{eqnarray}
\label{number}
  \reg(M/IM)=\max\{\reg(N/IN), \sat(IM)-1\}.
\end{eqnarray} By \ref{elem-1} we also have $\reg(N/IN) \leq \max\{\reg(IN)-1,
\reg(N)\}$, and since $N/IN$ is $0$-dimensional we conclude from the first part
of the proof  that $\reg(IN)\leq \reg(N)+\reg(I)$, so that
\begin{eqnarray*}
\reg(N/IN)&\leq &\max\{\reg(N)+\reg(I)-1, \reg(N)\}\\
&=&\reg(N)+\reg(I)-1\leq \reg(M)+\reg(I)-1.
\end{eqnarray*} The last inequality holds since $x$ is almost regular on
$M$. Thus together with (\ref{number}) we obtain
\begin{eqnarray}
\label{number2}
\reg(M/IM)\leq \max\{\reg(M)+\reg(I)-1, \sat(IM)-1\}.
\end{eqnarray} Notice further that
\begin{eqnarray}
\label{number3}
\reg(IM)\leq \max\{\reg(M), \reg(M/IM)+1\}.
\end{eqnarray} We may assume that $\reg(IM)> \reg(M)$, because otherwise
nothing is to prove. But then (\ref{number3}) implies that $\reg(IM)\leq \reg(M/IM)+1$. Hence together with 
(\ref{number2}) we get
\[
\reg(IM)\leq \max\{\reg(M)+\reg(I),\sat(IM)\}.
\] 
The desired inequality follows from (\ref{satIM}). 
\end{pf}

\section{Regularity of products of ideals of linear forms}

The goal of this section is to prove the following:
\begin{theorem}
\label{pils} Let $R$ be a polynomial ring and let $I_1, I_2, \dots, I_d$
be non-zero ideals of $R$ generated by linear forms. Then the product
$I_1I_2\cdots I_d$ has a linear resolution, i.e. $$\reg(I_1I_2\cdots
I_d)=d.$$
\end{theorem}

To prove the theorem we need some preliminary results. Let us fix some
notation. For a subset $A$ of
$\{1,\dots,d\}$ we will set  $I_A=\sum_{j\in A} I_j$ and denote by
$|A|$ the cardinality of $A$.  We have:

\begin{lemma}
\label{primdec} Let $I_1, I_2, \dots, I_d$ be non-zero ideals of $R$
generated by linear forms. Then $$I_1\cdots I_d= \sect_A I_A^{|A|}$$
    is a (possibly redundant) primary decomposition of $I_1\cdots I_d$.
Here the intersection is extended to all the non-empty subsets $A$ of
$\{1,\dots,d\}$.
\end{lemma}

As a corollary of \ref{primdec} we have:

\begin{corollary}
\label{satdeg} Let $I_1, I_2, \dots, I_d$ be non-zero ideals of $R$
generated by linear forms.
    Then $$\sat(I_1I_2\cdots I_d)\leq d.$$
\end{corollary}

\begin{pf} of \ref{satdeg}: Set $J=I_1I_2\cdots I_d$. By virtue of
\ref{primdec}
    $J^{sat}=\sect_A I_A^{|A|}$ where the intersection is extended to all
the non-empty subsets $A$ of $\{1,\dots,d\}$ such that $I_A\neq \mm$. It
follows that
$J=J^{sat}\sect \mm^d$ if $\sum I_i=\mm$ and  $J=J^{sat}$, otherwise.
This implies that $\sat(J)\leq d$.
\end{pf}

Now we prove \ref{primdec}:

\begin{pf}  The ideal $I_A^{|A|}$ is obviously $I_A$-primary and hence
it suffices to prove that
$I_1\cdots I_d= \sect_A I_A^{|A|}$. Set
$J=I_1I_2\cdots I_d$.   Let
$J_i$ be the product of the $I_j$ with $j\neq i$. By induction on
$d$, it is enough to show that:
$$J=J_1\sect \dots \sect J_d \sect (\sum_{i=1}^d  I_i)^d.$$ We prove
this equality by induction on $d$ and on $\dim R$.  The critical
inclusion is
$\supseteq$.  We may assume that
$\sum I_i=\mm$ (otherwise   all the ideals live in a smaller polynomial
ring). It is also harmless to assume that the residue field is
infinite.   Summing up, what we have to prove is that if   $f$ is an
element in $J_1\sect \dots \sect J_d$ of degree
$\geq d$ then $f\in J$.
     As $J_i$ is a product of
$(d-1)$  ideals of linear forms, by induction we know that  Corollary
\ref{satdeg} holds for $J_i$ and hence $\sat(J_i)\leq d-1$ for all $i$.
Let $x$ be a linear form which is a non-zerodivisor on $R/J_j^{sat}$ for
all the $J_j$ of positive dimension.  The ideals
$J+(x)/(x)$ of $R/(x)$ is the  product of ideals of linear forms
$I_i+(x)/(x)$. So, arguing  modulo $x$ and using induction on $\dim R$,
we see that
$f\in J+(x)$.  Write
$f=h+xf_1$, with $h\in J$. Replacing $f$ with $f-h$ we may  assume from
the really beginning that
$f=xf_1$. Since $f=xf_1\in J_i$ and $\sat(J_i)\leq d-1$, by the choice
of $x$ we may deduce that $f_1$ itself is  in  $J_i$ for all $i$.  Now
since the sum of the $I_i$ is
$\mm$ we may write
$x=\sum_i x_i$ with
$x_i \in I_i$. Then we have $f=xf_1=\sum_i x_if_1$ and each  $x_if_1
\in I_iJ_i=J$ so that $f \in J$.
\end{pf}

We are ready to prove \ref{pils}

\begin{pf} Set $J=I_1\dots I_d$.  Since $J$ is generated in degree
$d$ our   task is to prove that  $\reg(J)\leq d$.  We   prove it  by
induction on the dimension of $R$ and on $d$. The claim is trivial   if
$\dim R=1$.  If $\dim R/J=0$   then the assertion is also trivial.    We
may hence  assume that $\dim R/J>0$.   Let $x$ be a   linear form which
is a non-zerodivisor modulo $J^{sat}$. By \ref{sat}   we have that
$\reg(J)=\max\{ \reg(J+(x)), \sat(J) \}$. Note that
$\reg(J+(x))=1+\reg(R/J+(x))$. Since $\reg(R/J+(x))$ can be interpreted
as the regularity of
$R/J+(x)$ as an $R/(x)$-module and the ideal $J+(x)/(x)$  of $R/(x)$ is
a product of ideals of linear forms  we have  $\reg(R/J+(x))=d-1$. It
follows that
$\reg(I+(x))= d$. Since by \ref{satdeg}  $\sat(J)\leq d$,  we are done.
\end{pf}

The primary decomposition of \ref{primdec} is in general far from being
irredundant. For example we have:

\begin{proposition}
\label{lingen} Let  $V_1,\dots,V_d$  be a family of subspaces   of
$R_1$ which is linearly general, i.e.   one has
\[
\dim \sum_{i \in A} V_{i}=\min\{\dim R_1,
\sum_{i\in A} \dim  V_{i}
\}.
\] for all the non-empty subsets $A$ of $\{1,\dots,d\}$. Assume that
$\sum_{i=1}^dV_i=R_1$. Let $I_i$ be the ideal generated by $V_i$. Then
$$I_1\cdots I_d=I_1\sect \dots \sect I_d \sect \mm^d$$ is a primary
decomposition of
$I_1\cdots I_d$.
\end{proposition}

\begin{pf} We have to show that all the terms $I_A^{|A|}$ with
$1<|A|<d$ in the primary decomposition
\ref{primdec} are superfluous.  For such an $A$ we distinguish two
cases.  If
$\sum_{i\in A} \dim V_i\leq \dim R_1$ then by assumption
$\dim \sum_{i\in A} V_i= \sum_{i\in A} \dim V_i$ which implies that
$\sect_{i\in A} I_i=\Pi_{i\in A} I_i$.  Hence $I_A^{|A|}$ contains
$\sect_{i=1}^d I_i$ and it is therefore superfluous.  If instead
$\sum_{i\in A} \dim V_i>\dim R_1$, then by assumption $I_A=\mm$ and
hence $I_A^{|A|}\supset \mm^d$.
\end{pf}

    On the other hand there are cases where all the
$2^d-1$ ideals appearing in the primary decomposition \ref{primdec} are
essential.

\begin{example} Let $R=K[x_1,\dots,x_d,y]$ and consider
$I_i=(x_i,y)$. Set $J=I_1\cdots I_d$. It is not difficult to show that
for any subset $A\subseteq
\{1,\dots,d\}$ one has
$J:m=(y)+(x_i : i \in A)=I_A$ where $m=y^{|A|-1}\Pi_{i\not\in A}  x_i$.
Hence each $I_A$ is an associted prime of $J$. Therefore the primary
decomposition given in
\ref{primdec} is irredundant in this case.
\end{example}

\begin{question}  After \ref{pils} it is natural to ask whether
$$\reg(I_1I_2\cdots I_d)\leq \reg(I_1)+\reg(I_2)+\dots +\reg(I_d)$$
holds for ideals
$I_i$ generated by regular sequences.  By \ref{generic}, this is true if
each $I_i$ is generated by generic forms.
\end{question}

\section{Modules  with linear quotients}

We say that a graded $R$-module  $M$  has {\em linear quotients}  if
$M$ admits  a minimal system of generators $m_1,\dots, m_k$ such that
for every
$t=1,\dots,k$ one has that  $\langle m_1,\dots,m_{t-1} \rangle:_R m_t$
is an  ideal of $R$ generated by linear forms.

Examples of ideals with linear quotients are strongly stable and
    squarefree strongly stable ideals. Other important classes will be
considered in  the next sections.

\begin{lemma}  If  $M$ has linear quotients then
\[
\reg(M)=\max\{\deg m\:  m  \mbox{  is a minimal generator of } M\}.
\]  In particular, if all generators of $M$ have the same degree, then
$M$  has a linear  resolution over $R$.
\end{lemma}

\begin{pf} Let $m_1,\dots,m_k$ be as in the definition of module with
linear quotients. Set $M_t=\langle m_1,\dots,m_{t} \rangle$. We have an
exact sequence
$$0\to M_{t-1} \to M_t \to M_t/M_{t-1}\to 0$$ and $M_t/M_{t-1}$ is of
the form
$R/I[-\deg(m_t)]$ with $I$ an ideal  of $R$ generated by linear forms.
Since
$\reg(R/I)=0$ it follows that $\reg(M_t)\leq \max\{\reg(M_{t-1}),
\deg(m_t)\}$ and hence,  by induction, the assertion follows.
\end{pf}

\begin{example}
\label{lq} The ideal $J=(a^2b,abc,bcd,cd^2)$ of \ref{not}  has linear
quotients, the successive colons   being:
$$(0), \quad   (a), \quad  (a),\quad  (b).$$ On the other hand there are
ideals with linear resolution and without linear quotients. The easiest
example is the ideal $I$ of
$2$-minors  of the matrix
$$
\left(
\begin{array}{ccc} a & b & c \\ b & c & d
\end{array}
\right ).
$$
$I$ has a linear resolution but it cannot have linear quotients since it
is a prime ideal and hence $(f):(g)=(f)$ for each $f\in I$ with
$\deg(f)=2$.
\end{example}

Note that for a monomial ideal $I$ to have  linear quotients (with
respect to the monomial generators) is a purely combinatorial property
and hence does not depend on the characteristic of the base field. On
the other hand  the minimal free resolution of a monomial ideal, and
hence its linearity,  depends, in  general, on the characteristic of the
base field.  This shows that also for monomial  ideals  to have linear
quotients is a stronger property than to have  a linear  resolution. The
(famous) example of the Stanley-Reisner ideal of a triangulation of the
real  projective plane (see for example \cite[pag.236]{BH}) gives an
example of square free  monomial ideal that, if the characteristic of
$K$ is not $2$, has  a linear resolution and does  not have linear
quotients.

We have seen that the property of having a linear resolution is not
preserved by taking  products or powers of  ideals. The same thing  can
happen for the property of having linear quotients:

\begin{example}
\label{nlin} We know from \ref{lq} that $J=(a^2b,abc,bcd,cd^2)$  has
linear quotients, but as we have seen  in \ref{not}, $(b,c)J$ does not
even have  a linear resolution. Also, the ideal $I=(a^2b, a^2c, ac^2,
bc^2, acd)$  has linear quotients, the quotients being
$$(0),  \quad (b),\quad   (a), \quad  (a), \quad  (c, a).$$ But the
minimal resolution of   $I^2$  begins with
$$   R^{24}(-7)\dirsum R(-8) \to R^{15}(-6)\to I^2\to 0 $$ and hence
$I^2$ cannot have linear quotients.
\end{example}

\begin{question}  We have seen that   a product of ideals of
linear forms   has a linear resolution. One may ask whether such
an ideal has even linear quotients.  In the next section we will see
that  this is the case  for  products of ideals of variables, see \ref{product1}.  For the
general case, we have tested   many examples with CoCoA, starting with  generic and with
special ideals of linear forms.  We have always found ideals with linear quotients.
\end{question}

\section{Polymatroidal ideals}

In this section we consider a class of monomial ideals with linear
quotients which  is closed  under the operation of taking products. The
theorems presented here correspond to analogue  theorems in matroid
theory.

Let $R=K[x_1,\ldots, x_n]$ be the polynomial ring. For a monomial ideal
$I\subset R$ we denote  by $G(I)$ the unique minimal set of monomial
generators, and for a monomial
$u=x_1^{a_1}\ldots  x_n^{a_n}$ we set $\nu_i(u)=a_i$ for $i=1,\ldots,n$.

\begin{definition}
\label{polymatroidal} A  monomial ideal $I\subset R$ is said to be {\em
polymatroidal} if all its generators have  the same degree and if it
satisfies the following exchange property:
\begin{itemize}
\item[]for all  $u,v\in G(I)$ and all $i$ with $\nu_i(u)>\nu_i(v)$,
there exists an integer
$j$ with $\nu_j(v)>\nu_j(u)$ such that $x_j(u/x_i)\in G(I)$.
\end{itemize}
\end{definition}

The name is explained by the fact that  the elements of $G(I)$
correspond to the basis of a  polymatroid, as defined in \cite{W}. If
$I$ is a squarefree ideal, then this set corresponds  to the basis of a
matroid. Hence squarefree polymatroidal ideals are also called {\em
matroidal}.

For the convenience of the reader we reproduce from \cite{HT} the proof
of the following important property of polymatroidal ideals.

\begin{proposition}
\label{reproduce} A polymatroidal ideal $I$ has linear quotients with
respect to the reverse lexicographical  order of the generators.
\end{proposition}

\begin{pf}  Let $u\in G(I)$, and let $J$ be the ideal generated  by all
$v\in G(I)$ with $v>u$ (in the reverse lexicographical order). Then
\[ J:u=(v/[v,u]\: v\in J ).
\] Thus in order to prove that $J:u$ is generated by monomials  of
degree $1$,  we have to show that for each $v>u$ there exists $x_j\in
J:u$ such that $x_j$ divides
$v/[v,u]$.

In fact, let
$u=x_1^{a_1}\cdots x_n^{a_n}$ and
$v=x_1^{b_1}\cdots x_n^{b_n}$. Since $v>u$,  there exists an integer
$i$ with $a_i > b_i$ and $a_k=b_k$ for $k=i+1,\ldots,n$,   and hence an
integer $j$ with $b_j>a_j$ such that
$u'=x_j(u/x_{i})\in I$.  Since $j < i$, we see that $u'\in J$, and from
the equation
$x_iu'=x_ju$ we deduce that $x_j\in J:u$. Finally,  since
$\nu_j(v/[u,v])=b_j-\min\{b_j,a_j\}=b_j-a_j>0$,  we have that $x_j$
divides
$v/[v,u]$.
\end{pf}

Though products of ideals with linear quotients need not to have linear
quotients,  we nevertheless have

\begin{theorem}
\label{product} Let $I$ and $J$ be polymatroidal monomial ideals. Then
$IJ$ is polymatroidal.
\end{theorem}

\begin{pf} Let $u$ and $v$ be two monomials of same degree. We set
\[ d(u,v)=\frac{1}{2}\sum_i|\nu_i(u)-\nu_i(v)|.
\] Note that this is an integer. We call $d(u,v)$ the distance between
$u$ and $v$. This function  satisfies the usual rules of a distance
function. In particular, one has
$d(u,v)=0$ if and  only if $u=v$.

Now let $u_1, u\in G(I)$ and $v_1,v\in G(J)$ and suppose that
$\nu_i(u_1v_1)>\nu_i(uv)$. Then  we may assume that
$\nu_i(u_1)>\nu_i(u)$. Hence there exists an integer $j_1$ such that
$\nu_{j_1}(u)>\nu_{j_1}(u_1)$ and $u_2=x_{j_1}(u_1/x_i)\in G(I)$.
Moreover we have
$d(u_2,u)<d(u_1,u)$.

If $\nu(v)\geq \nu(v_1)$ we are done, because then
$\nu_{j_1}(uv)>\nu_{j_1}(u_1v_1)$, and
\[ x_{j_1}(u_1v_1/x_i)=u_2v_1\in G(IJ).
\]  Otherwise $\nu_{j_1}(v_1)>\nu_{j_1}(v)$. Hence there exists $k_1$
with
$\nu_{k_1}(v)>\nu_{k_1}(v_1)$ and such that $v_2=x_{k_1}(v_1/x_{j_1})\in
G(J)$. Moreover we  have $d(v_2,v)<d(v_1, v)$.

If $\nu_{k_1}(u)\geq \nu_{k_1}(u_2)$, then
$\nu_{k_1}(uv)>\nu_{k_1}(u_2v_1)=\nu_{k_1}(x_{j_1}(u_1v_1/x_i))$.  Thus
if $k_1\neq i$, then $\nu_{k_1}(uv)>\nu_{k_1}(u_1v_1)$, and we are done
since
\[ x_{k_1}(u_1v_1/x_i)=u_2v_2\in G(IJ).
\] On the other hand, if $k_1=i$, then $u_1v_1=u_2v_2$, and by induction
we may assume that the  exchange property holds since
$d(u_2,u)<d(u_1,u)$ and
$d(v_2,v)<d(v_1,v)$.

Otherwise $\nu_{k_1}(u_2)>\nu_{k_1}(u)$. Hence there exists $j_2$ with
$\nu_{j_2}(u)>\nu_{j_2}(u_2)$ and such that
$u_3=x_{j_2}(u_2/x_{k_1})\in G(I)$.  If
$\nu_{j_2}(v)\geq \nu_{j_2}(v_2)$, then
$\nu_{j_2}(uv)>\nu_{j_2}(u_2v_2)=\nu_{j_2}(x_{k_1}(u_1v_1/x_i))$. Thus
if $j_2\neq i$, then
$\nu_{j_2}(uv)>\nu_{j_2}(u_1v_1)$, and we are done since
\[ x_{j_2}(u_1v_1/x_i)=u_3v_2\in IJ.
\] On the other hand, if $j_2=i$, then $u_3v_2=u_1v_1$, and by induction
on the distance we have  the desired exchange property. Otherwise
$\nu_{j_2}(v_2)>\nu_{j_2}(v)$.

We may proceed in this way. Suppose we have already constructed sequences
$x_{j_1},\ldots,  x_{j_r}$,  $x_{k_1},\cdots, x_{k_{r-1}}$, and
$u_1,\ldots, u_{r+1}\in G(I)$, $v_1,\ldots,  v_r\in G(J)$ such that for
$i=1,\ldots,r$ we have
\begin{enumerate}
\item[(i)] $x_{k_{i-1}}$ divides $u_i$ and $x_{j_i}$ divides $v_i$,
\item[(ii)] $u_{i+1}=x_{j_i}(u_i/x_{k_{i-1}})$ and
$v_i=x_{k_{i-1}}(v_{i-1}/x_{j_{i-1}})$,
\item[(iii)] $d(u_{i+1},u)<d(u_i,u)$ and for $i\neq r$,
$d(v_{i+1},v)<d(v_i,v)$,
\item[(iv)] $\nu_{j_i}(u)>\nu_{j_i}(u_i)$ and
$\nu_{k_i}(v)>\nu_{k_i}(v)$.
\end{enumerate} Here we have set $k_0=i$ for systematic reasons. Notice
that
\[ u_{r+1}=x_{j_r}\cdots x_{j_1}(u_1/x_ix_{k_1}\cdots
x_{k_{r-1}})\quad\text{and}\quad v_r=x_{k_{r-1}}\cdots
x_{k_1}(v_1/x_{j_{r-1}}\cdots x_{j_1}),
\]

If $\nu_{j_r}(v)\geq \nu_{j_r}(v_r)$, then by (iv),
$\nu_{j_r}(uv)>\nu_{j_r}(u_rv_r)=\nu_{j_r}(x_{k_{r-1}}(u_1v_1/x_i))$.
Thus, if
$j_r\neq i$,  then $\nu_{j_r}(uv)>\nu_{j_r}(u_1v_1)$, and we are done
since
\[ x_{j_r}(u_1v_1/x_i)= u_{r+1}v_r\in G(IJ).
\] On the other hand, if $j_r=i$, and then $u_1v_1=u_{r+1}v_r$ and by
induction on the distance  we have the desired exchange property.

Otherwise $\nu_{j_r}(v_r)>\nu_{j_r}(v)$, and there exists $k_r$ with
$\nu_{k_r}(v)>\nu_{k_r}(v_r)$ and such that
$v_{r+1}=x_{k_r}(v_r/x_{j_r})\in G(J)$. Moreover  we have
$d(v_{r+1},v)<d(v_r,v)$. Thus the new elements $x_{k_r}$ and
$v_{r+1}$ satisfy again the properties (i)-(iv).

If $\nu_{k_r}(u)\geq \nu_{k_r}(u_{r+1})$, then by (iv),
$\nu_{k_r}(uv)>\nu_{k_r}(u_{r+1}v_r)=\nu_{k_r}(x_{j_r}(u_1v_1/x_i))$.
Thus, if
$k_r\neq i$, then $\nu_{k_r}(uv)>\nu_{k_r}(u_1v_1)$, and we are done
since
\[ x_{k_r}(u_1v_1/x_i)= u_{r+1}v_{r+1}\in G(IJ).
\] On the other hand, if $k_r=i$, and then $u_1v_1=u_{r+1}v_{r+1}$ and
by induction on the  distance we have the desired exchange property.

Otherwise $\nu_{k_r}(u_{r+1})>\nu_{k_r}(u)$, and there exists $j_{r+1}$
with
$\nu_{j_r}(u)>\nu_{j_r}(u_{r+1})$ and such
that$u_{r+2}=j_{r+1}(u_{r+1}/x_{k_r})\in G(I)$.  Moreover,
$d(u_{r+1},u)<d(u_r,u)$. Thus we have the conditions (i)-(iv) as before
but $r$ replaced by $r+1$. Condition (iii) implies that the process must
terminate. This proves the  theorem.
\end{pf}

Since ideals generated by subsets of the variables are obviously
polymatroidal, Theorem~\ref{product} implies

\begin{corollary}
\label{product1} Let $I_1,\ldots, I_d$ be ideals generated by subsets of
the variables. Then  $I=I_1\cdots I_d$  has linear quotients.
\end{corollary}

Let $I$ and $J$ be matroidal ideals. We let $I*J$ be the ideal which is
generated by all  monomials $uv$ with $u\in G(I)$ and $v\in G(J)$ such
that $uv$ is squarefree. We call $I*J$  the squarefree product of $I$
and $J$. Analogously to
\ref{product} we have

\begin{theorem}
\label{analogue} Let $I$ and $J$ be matroidal ideals. Then $I*J$ is
matroidal.
\end{theorem}

The proof of this theorem  similar to that of \ref{product}. We leave it
to the reader.

As a particular case of \ref{analogue} one has that the squarefree
product of ideals generated  by variables is matroidal. The
corresponding matroid is usually called {\em transversal}.

\section{Products of ideals defined by Hankel matrix}

In this section we use the notion  of ideals with linear quotients to
show that   products of   ideals of minors of  a Hankel matrix  have a
linear resolution.

Let $S$ be the polynomial ring $K[x_1,\dots,x_n]$ over some field
$K$.  Let  $X$  be a Hankel matrix with distinct entries
$x_1,\dots,x_n$; this means that  $X$ is an $a\times b$ matrix
$(x_{ij})$  with $x_{ij}=x_{i+j-1}$ and $a+b-1=n$. Let $I_t$ be the
ideal generated by the minors of size $t$ of $X$. It is known that $I_t$
does not depend on the size of the matrix $X$ (provided, of course,  $X$
contains $t$-minors); it depends only on
$t$ and $n$. For a given $n$ it follows that $t$ may vary from
$1$ to $m$, where $m=[(n+1)/2]$  is  the integer part of
$(n+1)/2$.   It is known that the powers of $I_2$ have a linear
resolution,
\cite{C2}.    Blum \cite[3.6]{B} has recently shown that if the Rees
algebra $R(I)$ of an ideal $I$ is Koszul then all the powers of $I$ have
linear resolutions. As we know that $R(I_t)$ is Koszul \cite{C2}, we
have that  $I_t^k$ has  a linear resolution for all $t$ and $k$.  We
prove here a  stronger result:

\begin{theorem} Let $X$ be a generic Hankel matrix.   Let
$t_1,\dots,t_p$ be integers and $I$ be the product of $I_{t_1}\cdots
I_{t_p}$. Then $I$ has a linear resolution.
\end{theorem}

We   recall some definitions and results from \cite{C}. Let $\tau$ be
the lexicographic term order on the monomials of $S$ and  $>_1$ the
partial order on
$x_1,\dots,x_n$ defined by $x_j>_1 x_i$ if and only if
$j-i>1$. A $>_1$-chain is a monomial  $x_{i_1}\cdots x_{i_k}$ such that
$x_{i_1}<_1 x_{i_2} <_1 \dots  <_1 x_{i_k}$. Denote by $J$ the initial
ideal of
$I=I_{t_1}\cdots I_{t_p}$  and by   $J_k$ that of  $I_k$. We know that
$$J_k=( m :  m \mbox{ is a  $>_1$-chain of degree }  k\}$$
    and that
$$J=J_{t_1}\cdots J_{t_p}.$$ Since the regularity can only increase by
passing to the initial ideal, it suffices to show that

\begin{proposition}
\label{Hlinq}  The ideal  $J$ has linear quotients.
\end{proposition}

Before proving \ref{Hlinq} we will describe the generators of $J$.  They
have a description in terms of the $\gamma$-functions associated to the
canonical decomposition of any monomial of $S$. Let us recall how.  Any
monomial $m$ of $S$ has a canonical decomposition
$m=m_1\cdots m_k$ as a product of $>_1$-chains. The monomial $m_1$ is
defined to be the largest, with respect to $\tau$,   among all the
$>_1$-chains which divide $m$. Similarly,   $m_2$ is  the largest among
all the
$>_1$-chains which divide $m/m_1$ and so on.  The shape of a monomial
$m$ is the sequence of integers $s(m)=\deg(m_1),\dots, \deg(m_k)$ where
$m=m_1\cdots m_k$  is the canonical decomposition of $m$. By the very
definition, the shape of
$m$ is a weakly decreasing sequence.  For any
$t$ and for any sequence of integers $s=s_1,\dots,s_p$ one defines
$$\gamma_t(s)=\sum_{i=1}^p \max(s_i-t+1,0).$$
    Furthermore, if $m$ is a monomial then we set:
    $$\gamma_t(m) =\gamma_t(s(m)). $$

\begin{example} Let $m=x_1^2x_2^3x_3^2x_5^3x_6x_7x_8^3$. Then
$$m= (x_1x_3x_5x_7)(x_1x_3x_5x_8)(x_2x_5x_8)(x_2x_6x_8)(x_2)$$ is the
canonical decomposition of $m$. Its shape is $s(m)=4,4,3,3,1$ and its
$\gamma$-values are $\gamma_1(m)=15,\gamma_2(m)=10, \gamma_3(m)=6,
\gamma_4(m)=2, and \gamma_t(m)=0$  for $t>4$.
\end{example}

Given the numbers $t_1,\dots, t_p$, let us denote by $\Omega$ the set of
the monomials
$m$ such that  $\deg(m)=\sum_{j=1}^p  t_j$ and
$\gamma_i(m)\geq \gamma_i(t_1,\dots,t_p)$ for every $i$.  In \cite{C} it
is proved:

\begin{proposition}
\label{hankel1}
\begin{itemize}
\item[(1)]  $\Omega$ is a system of generators of $J$,
\item[(2)]  Let  $m$ be a monomial with a decomposition (canonical or
not)  $m=n_1\cdots n_v$ where the
$n_i$ are $>_1$-chains.  Set   $s=\deg(n_1), \dots, \deg(n_v)$. Then
$\gamma_i(m)\geq
\gamma_i(s)$  for every  $i$.
\end{itemize}
\end{proposition}

We introduce a total order  $\sigma$ on the monomials of $S$  as
follows. Let  $m,n$ be monomials  of  $S$ and $m=m_1\cdots m_k$ and
$n=n_1\cdots n_h$ their canonical decompositions. We set $m>_\sigma n$
if $m_j>_\tau n_j$   for  the first index $j$ such that $m_j\neq n_j$.
Note that
$\sigma$ is different from $\tau$; for instance $x_1^2>_\tau x_1x_3$ but
$x_1x_3>_\sigma x_1^2$. Note also that $\sigma$ is not a term order. Now
we are ready to prove:

\begin{pf} of \ref{Hlinq}: We show that $J$ has linear quotients with
respect to the set of generators $\Omega$   totally ordered by
$\sigma$. Let $m,n$ be elements of $\Omega$ with $m>_\sigma n$. We have
to show that there exists $v\in \Omega$ such that $v>_\sigma n$,
$v/[v,n]$ divides $m/[m,n]$  and $\deg [v,n]= \deg v-1$. Let
$m=m_1\cdots m_k$ and  $n=n_1\cdots n_h$ be the  canonical
decompositions and let $j$ be the smallest index such that $m_j\neq
n_j$. Then $m_j>_\tau n_j$. Let
$m_j=x_{a_1}\cdots x_{a_r}$ and  $n_j=x_{b_1}\cdots x_{b_s}$. Then there
exists a index
$z$ such that   $a_i=b_i$ for $i=1,\dots, z-1$ and either $a_z<b_z$ or
$s=z-1$ and
$r\geq z$. In the  former  case ($a_z<b_z$)  we put
$v=nx_{a_z}/x_{b_z}$. In the latter case  we put
$v=nx_{a_z}/x_{q}$ where $x_q$ is a variable which appear in
$n_{j+1}$ (note that
$h>j$, since $m$ and $n$ have both degree $\sum t_i$). We have  to show
that $v$ has the desired properties.

First of all,  note that
$v/[v,n]=x_{a_z}$.  This is clear in the first case while in the second
it follows from the fact that $q$ cannot be equal to  $a_z$ otherwise the
$j$-th factor in the canonical decomposition of $n$ would be  a multiple
of
$x_{b_1}\cdots x_{b_{z-1}}x_{a_z}$.

Secondly,  we claim that   $x_{a_z}$ divides $m/[m,n]$. To this end,
note that
$m/[m,n]=m^{'}/[m^{'},n^{'}]$ where
$m^{'}=m/e$ and  $n^{'}=n/e$ and $e$ is the common initial part of the
canonical decomposition, i.e. $e=m_1\cdots m_{j-1}x_{a_1}\cdots
x_{a_{z-1}}$. Since $x_{a_z}$ appears in  $m^{'}$ and it does not appear
in
$n^{'}$ (otherwise, as above,  the $j$-th factor in the canonical
decomposition of
$n$ would be  a multiple of   $x_{b_1}\cdots x_{b_{z-1}}x_{a_z}$), we
may conclude that  $x_{a_z}$ divides
$m/[m,n]$.

It remains to show that $v$ belongs to $\Omega$ and that $v>_\sigma n$.
In the case
$a_z<b_z$ note that the $v$ has a decomposition into $>_1$-chains
$v=n_1\cdots n_{j-1}un_{j+1}\cdots n_h$ with $u=n_ix_{a_z}/x_{b_z}$.
This need not to be the canonical decomposition, but its  shape is equal
to that  of the canonical decomposition of $n$ and this is enough (by
\ref{hankel1}) to conclude that $v\in
\Omega$. Since by construction  $u>_\tau n_j$, it is not difficult to
check that
$v>_\sigma n$.  In the case $s=z-1$ and $r\geq z$ note that the $v$ has
a decomposition into $>_1$-chains $v=n_1\cdots
n_{j-1}u_1u_2n_{j+2}\cdots n_h$ with
$u_1=n_jx_{a_z}$ and
$u_2=n_{j+1}/x_{a_z}$. As above, this need not to be the canonical
decomposition.  Its shape has been obtained from the shape of $n$ by
the operation ``increase a larger factor and decrease a shorter". The
effect of this operation on the
$\gamma$-values is clear: the $\gamma$-values cannot decrease. This,
together with the fact that $n$ is in $\Omega$ and
\ref{hankel1} implies that $v$ is in $\Omega$.  As in the other case,
since
$u_1>_\tau n_j$ one can also deduce that $v>n$.
\end{pf}

\end{document}